\documentclass[preprint,review,12pt]{elsarticle}

\usepackage{amssymb}

\journal{journal}
\begin{document}

\begin{frontmatter}
\title{An analytical proof for Lehmer's totient conjecture using Mertens' theorems}
\author{Ahmad Sabihi} 
\ead{sabihi2000@yahoo.com}
\address{Teaching professor and researcher at some universities of Iran}

\begin{abstract}
We make an analytical proof for Lehmer's totient conjecture. Lehmer conjectured that there is no solution for the congruence equation $n-1\equiv 0~(mod~\phi(n))$ with composite integers,$n$, where $\phi(n)$ denotes Euler's totient function. He also showed that if the equation has any composite solutions, $n$ must be odd, square-free, and divisible by at least 7 primes. Several people have obtained conditions on values ,$n$, and number of square-free primes constructing $n$ if the equation can have composite solutions. Using Mertens' theorems, we show that it is impossible that the equation can have any composite solution and implies that the conjecture should be true for all the positively composite numbers.           
\end{abstract}

\end{frontmatter}

\section{Introduction}
 Lehmer's totient conjecture was stated by D.H. Lehmer in 1932 \cite{L}. Lehmer conjectured that there are no composite solutions,$n$, for the equation $n-1\equiv 0~(mod~\phi(n))$ . We know that this conjecture is true for every prime numbers. He also proved that if any such ,$n$, exists, it must be odd, square-free, and divisible by at least seven primes \cite{L}. Pinch calls such an $n$ a \textit{Lehmer number} and defines the \textit{Lehmer index} of $n$ to be the ratio $\frac{n-1}{\phi(n)}$ \cite{P}. As we should know every exponent $\lambda(n)$ of the multiplication group $(\mathbb Z/\mathbb N )^*$ must divide $n-1$ and follows that a Carmichael number $n$ must be square-free with at least three prime factors, and $p-1|n-1$ for every prime $p$ dividing $n$. Conversely, any such $n$ must be a Carmichael number. Since the exponent $\lambda(n)$ of the multiplicative group divides its order $\phi(n)$, a Lehmer number must be a Carmichael number. Lieuwens \cite{Li} showed that a Lehmer number divisible by 3 must have index at least 4 and hence must have at least 212 prime factors and exceeds $5.10^{570}$. Kishore \cite{K} proved that a Lehmer number of index at least 3 must have at least 33 prime factors and exceeds $2.10^{56}$. Cohen and Hagis \cite{CH} showed that a Lehmer number divisible by 5 and of index 2 must have at least 13 prime factors and if we have any composite solution ,$n$, to the problem, then $n>10^{20}$ and number of prime factors must be greater than or equal 14. We firstly show that using Mertens' theorems, we are able to asymptotically prove that the equation $n-m\phi(n)=1$ with odd composite number ,$n$, and having k square-free prime factors cannot have any solutions. We also investigate about the equation  $n-m\phi(n)=-1$ and take a conclusion that this equation may have solutions as Lehmer has shown in his paper \cite{L}.
We decompose our proof into the four theorems 3 to 6. Then, we show that ,$n$, must be odd, and square-free as Lehmer showed before, but by another method. To prove our theorems, we make use of Mertens' theorems on the density of primes and re-prove some of them.  
\setcounter{equation}{0}
\section{Theorems} \label{sc2}
\subsection{\textbf{Theorem 1: Mertens' 2nd theorem \cite{M}}}
\textit{Let $p$ be a prime and $x>1$ every real number, then 
\begin{equation}
\sum_{p\leq x} \frac{1}{p}=\log\log x+a+O(\frac{1}{\log x})
\end{equation}
where a possible value of "$a$" can be $a=0.2614972128...$}
\subsection{\textbf{Theorem 2: Mertens' 3rd theorem \cite{M}}}
\textit{Let $p$ be a prime and $x>1$ every real number, then 
\begin{equation}
\prod_{p\leq x}(1-\frac{1}{p})\sim\frac{e^{-\gamma}}{\log x}
\end{equation}
where the notation $f(x)\sim g(x)$ means that limitation $\frac{f(x)}{g(x)}=1$ when x tends to infinity. $\gamma$ denotes Euler's constant.} 
\subsection{\textbf{Corollary 1:}}
\textit{Let $p$ be a prime, $x>1$ every real number, and $c>0$ an absolute constant, then 
\begin{equation}
\prod_{p\leq x}(1-\frac{1}{p})>\frac{c}{\log x}
\end{equation}
where "$c$" can be 0.3 for $x\geq 2973$ and 0.09 for $x\geq 3$ in this paper.}
\subsection{\textbf{Theorem 3:}}
\textit{Let $p_{i}$ to $p_{k}$ be all of the prime factors including only odd square-free prime factors of the odd number $n$ and sufficiently so large integers or all of prime factors values tend to infinity versus the number of them, then the equation
\begin{equation}
n-m\prod_{p_{i}\leq p\leq p_{k}}(p-1)=\pm1
\end{equation}
does not any solution. $m$ denotes a positive integer.}
\subsection{\textbf{Theorem 4:}}
\textit{Let $p_{1}$ to $p_{k}$ be all of the prime factors including only odd square-free prime factors of the odd number $n$,all of them be existed, and  $p_{k}$  sufficiently so large integer or tends to infinity, then the equation
\begin{equation}
n-m\prod_{p\leq p_{k}}(p-1)=\pm1
\end{equation}
does not any solution.  $m$ denotes a positive integer.}
\subsection{\textbf{Theorem 5:}}
\textit{Let $p_{i}$ to $p_{k}$ be all of the prime factors including only odd square-free prime factors of the odd number $n$ and $p_{k}$ sufficiently so large integer or tends to infinity, then 
\begin{equation}
n-m\prod_{p_{i}\leq p\leq p_{k}}(p-1)=\pm1
\end{equation}
does not any solution. $m$ denotes a positive integer.}
\subsection{\textbf{Theorem 6:}}
\textit{Let $p_{i}$ to $p_{k}$ be all of the prime factors including only odd square-free prime factors of the odd number ,$n$, and none of them be so large and unbounded (all of them be bounded), then the equation    
\begin{equation}
n-m\prod_{p_{i}\leq p\leq p_{k}}(p-1)=1
\end{equation}
does not any solution, but the equation
\begin{equation}
n-m\prod_{p_{i}\leq p\leq p_{k}}(p-1)=-1
\end{equation}
may have solutions. $m$ denotes a positive integer.}
\setcounter{equation}{0}
\section{Proofs} \label{sc2}
\subsection{\textbf{Proof of Theorem 1}}
As is well-known, Mertens himself has proven this theorem but we give another method for making its proof. We really reprove (reformulate) the proof. 
The proof can be made by applying three times \textit{the Abel summation formula} to the the series
\begin{equation}
\sum_{p\leq x} \frac{1}{p}=\sum_{p\leq x} \frac{\log p}{p}.\frac{1}{\log p}
\end{equation}
Firstly, we apply it to the series $\sum_{p\leq x}\log p$ and reach to $\theta(x)=x+o(\frac{x}{\log x})$. Secondly, $\sum_{p\leq x} \frac{\log p}{p}$. Let $\sum_{p\leq x} \log p=\theta(x)$ and $\phi(x)=\frac{1}{x}$ and substitude them into the Abel summation formula as follows:
\begin{equation}
\sum_{p\leq x} \frac{\log p}{p}=\theta(x)\phi(x)+\int_{1}^{x}\frac{\theta(u)}{u^{2}}du
\end{equation}
Then, we have
\begin{equation}
\sum_{p\leq x} \frac{\log p}{p}=1+\log(x)+o(\log\log x)
\end{equation}
Thirdly, we apply it to the entire series.
 Let $A(x)=\sum_{p\leq x} \frac{\log p}{p}$ and $\Phi(x)=\frac{1}{\log(x)}$ into the Abel summation formula
\begin{eqnarray}
\sum_{p\leq x} \frac{1}{p}=A(x)\Phi(x)+\int_{2}^{x}A(u).\Phi'(u).du=\{1+\log x+o(\log\log x)\}.(\frac{1}{\log x})\nonumber\\+\int_{2}^{x}\frac{\{1+\log u+o(\log\log u)\}}{u(\log u)^{2}}du=1+\log x+o(\frac{\log\log x}{\log x})+\int_{2}^{x}\frac{du}{u(\log u)^{2}}+\nonumber\\ \int_{2}^{x}\frac{du}{u\log u}+\int_{2}^{x}o(\frac{\log\log u}{u(\log u)^{2}})=1+\frac{1}{\log 2}-\log\log 2+\log\log x+o(\frac{\log\log x}{\log x})+d+\nonumber\\ o(\frac{\log\log x}{\log x}+\frac{1}{\log x})=1+\frac{1}{\log 2}-\log\log 2+d+\log\log x+o(\frac{1}{\log x})~~~~~~
\end{eqnarray}
Where $d$ denotes all unknown constant values created in (3.4). Since according to the properties of small "o" and big "O" notations, we have $o(\frac{1}{\log x})=O(\frac{1}{\log x})$,then
\begin{equation}
\sum_{p\leq x} \frac{1}{p}=\log\log x+a+O(\frac{1}{\log x})
\end{equation}
where $a=1+\frac{1}{\log 2}-\log\log 2+d$. Although, precisely calculating $a$ is difficult, but our attempts to calculate the value $a$ using directly processing data by substituting into (3.5) gave us an approximate value about 0.261497...  
\subsection{\textbf{Proof of Theorem 2}}
The proof can be found in the Mertens' paper \cite{M}.
\subsection{\textbf{Proof of Corollary 1}}
The proof can easily be made by appealing to the Riemann Zeta Function and Euler's product \cite{WE} as follows:
\begin{equation}
\zeta^{-1}(s)=\prod_{p}(1-\frac{1}{p^{s}})=\sum_{n=1}^{\infty}\frac{\mu(n)}{n^{s}}
\end{equation}
Putting $s=1$ in (3.6), we have
\begin{equation}
\prod_{p}(1-\frac{1}{p})=\sum_{n=1}^{\infty}\frac{\mu(n)}{n}
\end{equation}
and trivially checking gives us
\begin{equation}
\prod_{p\leq x}(1-\frac{1}{p})>\prod_{p}(1-\frac{1}{p})=\sum_{n=1}^{\infty}\frac{\mu(n)}{n}
\end{equation}
Abel Summation Formula gives us again that assuming $\sum_{n\leq x}\mu(n)=o(x)$ \cite{WE}, we have
\begin{equation}
\sum_{n=1}^{\infty}\frac{\mu(n)}{n}=\frac{e^{-\gamma}}{\log x}+o(\log(x))
\end{equation}
Combining (3.9) with theorem 2 and (3.8) we find
\begin{equation}
\prod_{p\leq x}(1-\frac{1}{p})>\prod_{p}(1-\frac{1}{p})=\sum_{n=1}^{\infty}\frac{\mu(n)}{n}=\frac{e^{-\gamma}}{\log x}+o(\log(x))>\frac{c}{\log(x)}
\end{equation}
If we let $c<e^{-\gamma}$, then inequality and the theorem is completed. We choose $c=0.3$ in this paper. On the other hand, we appeal to Theorem 7, Corollary of the Rosser and Schoenfeld's paper \cite{RS} (the relation (3.27)) and we find that for $x=3$, we can choose $c=0.09$  since the term $e^{-\gamma}(1-\frac{1}{(\log 3)^{2}})=0.0962709...$. Therefore, we choose a new lower bound for $c$ i.e. $c=0.09$ since we have $e^{-\gamma}(1-\frac{1}{(\log x)^{2}})>e^{-\gamma}(1-\frac{1}{(\log 3)^{2}})>0.09$ for all the odd primes $\geq 3$. Also, Dusart \cite{Du} in 2010, stated the Theorem 6.12 giving a new bound for all $x\geq 2973$. This new bound for $x\geq 2973$ is 0.46842432.... This means that $c=0.3$ is acceptable for these values as well.
\subsection{\textbf{Proof of Theorem 3}}
If we divide both of sides of the equation (2.4) by $\prod_{p_{i}\leq p\leq p_{k}}(p-1)$,then
\begin{equation}
\frac{1}{\prod_{p_{i}\leq p\leq p_{k}}(1-\frac{1}{p})}-m=\frac{\pm1}{\prod_{p_{i}\leq p\leq p_{k}}(p-1)}
\end{equation}
Since $n=p_{i}...p_{k}$ is odd and $p_{i}$ to $p_{k}$ are also odd square-free prime factors of $n$,then trivially all of them must be $\geq 3$ and follows $\frac{1}{\prod_{p_{i}\leq p\leq p_{k}}(p-1)}<\frac{1}{8}$ if the numerator of right side be (+1) and $-\frac{1}{\prod_{p_{i}\leq p\leq p_{k}}(p-1)}>-\frac{1}{8}$ if the numerator of right side be (-1). On the other hand, the left side of (3.11) should be greater than zero for the plus sign and less than zero for the minus sign. Therefore, for plus sign we have
\begin{equation}
\frac{1}{\prod_{p_{i}\leq p\leq p_{k}}(1-\frac{1}{p})}-\frac{1}{8}<m<\frac{1}{\prod_{p_{i}\leq p\leq p_{k}}(1-\frac{1}{p})}
\end{equation}    
and for minus sign
\begin{equation}
\frac{1}{\prod_{p_{i}\leq p\leq p_{k}}(1-\frac{1}{p})}<m<\frac{1}{\prod_{p_{i}\leq p\leq p_{k}}(1-\frac{1}{p})}+\frac{1}{8}
\end{equation}    
Since our assumption says us that all $p_{i}$ to $p_{k}$ tend to infinity versus the number of primes within the interval $(p_{i},p_{k})$, the relations (3.12) and (3.13) change to 
\begin{equation}
1+\varepsilon-\frac{1}{8}<m<1+\varepsilon
\end{equation}    
\begin{equation}
1+\varepsilon<m<1+\varepsilon+\frac{1}{8}
\end{equation} 
when $\varepsilon$ tends to zero.
This is due to if we let $M$ denotes the number of primes from $p_{i}$ to $p_{k}$ (note: there may not exist all of consecutive primes within the interval $(p_{i},p_{k})$) then we find 
\begin{equation}
(1-\frac{1}{p_{i}})^{M}\leq \prod_{p_{i}\leq p\leq p_{k}}(1-\frac{1}{p})\leq (1-\frac{1}{p_{k}})^{M}
\end{equation}
Since according to our assumption, $p_{i}$ to $p_{k}$ are so large versus $M$, then all the fractions $\frac{M}{p_{i}}$ to $\frac{M}{p_{k}}$ tend to zero and
\begin{equation}
\lim_{p_{i}\longrightarrow \infty}(1-\frac{1}{p_{i}})^{M}=\lim_{p_{i}\longrightarrow \infty}\{(1-\frac{1}{p_{i}})^{p_{i}}\}^{\frac{M}{p_{i}}}=\lim_{p_{i}\longrightarrow \infty}(\frac{1}{e})^{\frac{M}{p_{i}}}=1
\end{equation}
and in the similar way
\begin{equation}
\lim_{p_{k}\longrightarrow \infty}(1-\frac{1}{p_{k}})^{M}=1
\end{equation}
Then the inequality (3.16) gives us
\begin{equation}
\lim_{p_{i}~to~p_{k} \longrightarrow \infty}\prod_{p_{i}\leq p\leq p_{k}}(1-\frac{1}{p})=1
\end{equation} 
This means that the integer number $m$ can only be 1 when (3.14) holds and cannot be any integer number when (3.15) holds. If $m=1$, it is impossible to hold by appealing to Lehmer's paper \cite{L} since $m=1$ if and only if $n$ is prime. This completes the proof. 
\subsection{\textbf{Proof of Theorem 4}} 
If we divide both of sides of (2.5) by $n=p_{1}...p_{k}$ and substitute $\prod_{p}(1-\frac{1}{p})=\frac{e^{-\gamma}}{\log x}+o(\log(x))$ from (3.8) and (3.9) into it, then let $x=p_{k}$
\begin{equation}
1-m\{\frac{e^{-\gamma}}{\log p_{k}}+o(\log(x))\}=\frac{\pm1}{p_{1}...p_{k}}
\end{equation}
Since $p_{k}\longrightarrow \infty$ then also $x\longrightarrow \infty$ and (3.20) changes to
\begin{equation}
1-\frac{m}{e^{\gamma}\log p_{k}}=\frac{\pm1}{p_{1}...p_{k}}
\end{equation}
The right side tends to zero since $p_{k}\longrightarrow \infty$. This means that the left side should also tend to zero and $m$ is of order $e^{\gamma}\log x$. Since $m$ is a positive integer, it could be of the form
\begin{equation}
m=[e^{\gamma}\log p_{k}]=e^{\gamma}\log p_{k}-\alpha
\end{equation}
for plus sign since the left side of (3.21) should be closed to $0^{+}$  or
\begin{equation}
m=[e^{\gamma}\log p_{k}]+j=e^{\gamma}\log p_{k}+j-\alpha
\end{equation}
for minus sign since the left side of (3.21) should be closed to $0^{-}$ 
where $j\geq 1$ and denotes an integer, the sign $[~]$ denotes the integer part of a number, and $\alpha$ denotes the fractional part of $e^{\gamma}\log x$. Therefore, the relation (3.21) can be changed into
\begin{equation}
\frac{\alpha}{e^{\gamma}\log p_{k}}=\frac{1}{p_{1}...p_{k}}~~ or~~ \frac{j-\alpha}{e^{\gamma}\log p_{k}}=\frac{1}{p_{1}...p_{k}}
\end{equation}
for when $p_{k},x\longrightarrow \infty$.
Since the denominator of the right side fraction of (3.24) is of the order more than $p_{k}$ and the denominator of the left right fraction is of order $\log p_{k}$, $\alpha$ and $j-\alpha$ are also bounded, then these two sides cannot be equal for when $p_{k}$ is tending to infinity and the equation (3.24), (3.21), and finally (2.5) cannot have any solutions.
\subsection{\textbf{Proof of Theorem 5}} 
The proof of this theorem also likes to Theorem 4. Consider all primes $p_{i}$ to $p_{k}$ exist or missing some of them, then regarding Theorem 2
\begin{equation}
A(p).\prod_{p_{i}\leq p\leq p_{k}, p_{k}\longrightarrow \infty}(1-\frac{1}{p})=\frac{e^{-\gamma}}{\log p_{k}}
\end{equation}
where $A(p)$ denotes a function of prime numbers before $p_{k}$ or some before  $p_{k}$ and some between $p_{i}$ and $p_{k}$ for completing and converting $\prod_{p_{i}\leq p\leq p_{k}}(1-\frac{1}{p})$ to $\prod_{p\leq p_{k}}(1-\frac{1}{p})$,which may be a constant value or variative one.Similarly to (3.21), we have
\begin{equation}
1-\frac{m}{A(p)e^{\gamma}\log p_{k}}=\frac{\pm1}{p_{1}...p_{k}}
\end{equation}
where 
\begin{equation}
A(p)=\prod_{p\leq p_{(i-1)}}(1-\frac{1}{p})~~or~~A(p)=\prod_{p\leq p_{(i-1)}}(1-\frac{1}{p}).\prod_{p_{i}\leq  p_{m}\leq p_{k}}(1-\frac{1}{ p_{m}})
\end{equation}
and $p_{m}$ denotes primes missing within the interval $(p_{i}, p_{k})$. Trivially, $A(p)<1$. Similarly to the proof of Theorem 4, we find
\begin{equation}
m=[A(p)e^{\gamma}\log p_{k}]=A(p)e^{\gamma}\log p_{k}-\beta
\end{equation}
for plus sign since the left side of (3.26) should be closed to $0^{+}$  or
\begin{equation}
m=[A(p)e^{\gamma}\log p_{k}]+j=A(p)e^{\gamma}\log p_{k}+j-\beta
\end{equation}
since the left side should be closed to $0^{-}$ for when $p_{k}\longrightarrow \infty$. Therefore, (3.26) changes to 
\begin{equation}
\frac{\beta}{A(p)e^{\gamma}\log p_{k}}=\frac{1}{p_{1}...p_{k}}
\end{equation}
 or
\begin{equation}
\frac{j-\beta}{A(p)e^{\gamma}\log p_{k}}=\frac{1}{p_{1}...p_{k}}
\end{equation}
Where $0\leq \beta <1$. To being better closed to zero in relation (3.31), we should choose $j=1$. The arguments are similar to the arguments of Theorem 4 and the proof is completed.
\subsection{\textbf{Proof of Theorem 6}}
Regarding Corollary 1,we have
\begin{equation}
\frac{c}{\log p_{k}}<\prod_{p\leq p_{k}}(1-\frac{1}{p})=\prod_{p_{1}\leq p\leq p_{(i-1)}}(1-\frac{1}{p}).\prod_{p_{i}\leq  p_{m}\leq p_{k}}(1-\frac{1}{ p_{m}}).\prod_{p_{i}\leq p\leq p_{k}}(1-\frac{1}{p}) 
\end{equation}
where $p_{m}$ denotes primes missing within the interval $(p_{i}, p_{k})$.
Let $A(p)=\prod_{p_{1}\leq p\leq p_{(i-1)}}(1-\frac{1}{p}).\prod_{p_{i}\leq  p_{m}\leq p_{k}}(1-\frac{1}{ p_{m}})$ and knowing $A(p)<1$ then
\begin{equation}
\frac{c}{A(p)\log p_{k}}<\prod_{p_{i}\leq p\leq p_{k}}(1-\frac{1}{p}) 
\end{equation} 
and multiplying the left side by a coefficient $l_{k}>1$, we find an equation
\begin{equation}
\frac{cl_{k}}{A(p)\log p_{k}}=\prod_{p_{i}\leq p\leq p_{k}}(1-\frac{1}{p}) 
\end{equation}
Similarly to (3.26), we have
\begin{equation}
1-\frac{mcl_{k}}{A(p)\log p_{k}}=\frac{\pm1}{p_{i}...p_{k}}
\end{equation}
As the Theorems 4 and 5 arguments, we have
\begin{equation}
m=[A(p)\frac{\log p_{k}}{cl_{k}}]=\frac{A(p)\log p_{k}}{cl_{k}}-\psi
\end{equation}
for plus sign since the left side of (3.35) should be closed to $0^{+}$  or
\begin{equation}
m=[A(p)\frac{\log p_{k}}{cl_{k}}]+j=\frac{A(p)\log p_{k}}{cl_{k}}+j-\psi
\end{equation}
for minus sign to be closed to $0^{-}$ (According to Lehmer's, Cohen's, Kishore's, and Lieuwens' arguments, the number of prime numbers to have a composite solution should be more than 7,14,33, or 212. Thus,the value $\frac{1}{p_{i}...p_{k}}$ should be certainly closed to zero due to being large the value $p_{i}...p_{k}$). $0\leq\psi<1$ denotes the fractional part of a positive   real number.Therefore, the equations (2.7) and (2.8) are found respectively
\begin{equation}
\frac{\psi cl_{k}}{A(p)\log p_{k}}=\frac{1}{p_{i}...p_{k}}  
\end{equation}
and
\begin{equation}
\frac{(j-\psi) cl_{k}}{A(p)\log p_{k}}=\frac{1}{p_{i}...p_{k}}  
\end{equation}
As we know, $\log p_{k}$ isn't an integer number and since $cl_{k}>0.09$ regarding Corollary 1  and $A(p)$ tends to zero by increasing the number of primes and being larger $p_{(i-1)}$ and $p_{k}$, then the value  $\frac{A(p)}{cl_{k}}$ gets smaller and smaller and $\log p_{k}$ larger and larger, thus the fractional part of $A(p)\frac{\log p_{k}}{cl_{k}}$ gets closer to the number 1. This means that $\psi$ gets closer to the number 1 to that of zero. Therefore, since $\frac{\psi cl_{k}}{A(p)}$ gets larger than 1,then (3.38) cannot have any solution, but may (3.39) have solution since $(j-\psi)$ gets closer to zero with $j=1$ and $\frac{(1-\psi) cl_{k}}{A(p)}$ gets closer to zero as well. As Lehmer \cite{L}, Kishore \cite{K}, Cohen \cite{CH}, and specifically Lieuwens \cite{Li} showed that if the case $n-m\prod_{p_{i}\leq p\leq p_{k}}(p-1)=1$ has composite solution,then the number of prime factors should be at least 7, 14, 33 or 212, therefore, we see that the order of magnitude of $n$ must be very large and our hypothesis can be more precise.
 
\textbf{Example:}

Lehmer showed that $n=3.5.17.257$ is a composite solution for the equation $n-m\prod_{p_{i}\leq p\leq p_{k}}(p-1)=-1$. We compute the values $1-\psi$, $A(p)$, $cl_{k}$, and $\log p_{k}$ assuming $c=0.09$ and substitute them into (3.38) and (3.39) as follows:\\
Here, we have $p_{i}=3$, $p_{i+1}=5$ ,  $p_{i+5}=17$,  $p_{k}=257$. For computing $A(p)$, one should compute all of other missing primes as:
\begin{eqnarray}
A(p)=(1-\frac{1}{7})(1-\frac{1}{11})(1-\frac{1}{13})(1-\frac{1}{19})(1-\frac{1}{23})...\nonumber\\ (1-\frac{1}{211})(1-\frac{1}{223})(1-\frac{1}{227})(1-\frac{1}{229})(1-\frac{1}{233})~~~~\nonumber\\(1-\frac{1}{239})(1-\frac{1}{241})(1-\frac{1}{251})=0.39984516~~~~~
\end{eqnarray}   
The relation (3.34) gives us $\prod_{p_{i}\leq p\leq p_{k}}(1-\frac{1}{p})=0.50000763$,$\log p_{k}=\log 257=5.54907608$,$l_{k}=12.3266964$, and $\psi=0.99996922$. Just, we are ready to compute the left sides of the two relations (3.38) and (3.39).

From the left side of the relation (3.38), we find
\begin{equation}
\frac{\psi cl_{k}}{A(p)\log p_{k}}=0.49999238 
\end{equation}
and from the left side of (3.39) with $j=1$ we find
\begin{equation}
\frac{(1-\psi) cl_{k}}{A(p)\log p_{k}}=1.526023\times 10^{-5}
\end{equation}
If we compute the right side of each of two relations (3.38) and (3.39)
\begin{equation}
\frac{1}{p_{i}...p_{k}}=\frac{1}{3\times 5\times 17\times 257}=1.525902\times 10^{-5}  
\end{equation}
and compare to the corresponding left sides (the relations (3.41) and (3.42)), then we find that the equation (2.8) has a solution since the left and right sides are very close to each other, but the equation (2.7) (same Lehmer's conjecture) does not any solution since the left and right sides are far from each other. Another example can be made by other composite number $n=3.5.17.257.65537$, which Lehmer showed it can be a composite solution for the same equation. Since $p_{k}=65537\geq 2973$, we consider $c=0.3$ for our calculations.       
\subsection{\textbf{Lehmer's totient conjecture}}
We discuss about Lehmer's totient conjecture here. Firstly, we know if $n$ is a prime number ,$p$, then $\phi(n)=p-1$ and trivially implies $n-1\equiv 0~(mod~\phi(n))$. Conversely, if we have  $n-1\equiv 0~(mod~\phi(n))$,then let $n=p_{1}^{t_{1}}...p_{k}^{t_{k}}$ be prime factors decomposition of $n$. This means that $\phi(n)=p_{1}^{t_{1}}...p_{k}^{t_{k}}\prod_{p_{1}\leq p\leq p_{k}}(1-\frac{1}{p})$. If $t_{1},...,t_{k}\geq 2$, then we find that $p_{1},...,p_{k}| both~\phi(n)~~and~~ n$. On the other hand, regarding our assumption $n-1\equiv 0~(mod~\phi(n))$ and we should have $p_{1},...,p_{k}| (n-1)$. But, these imply that $p_{1},...,p_{k}| gcd(n,n-1)=1$, which is impossible to occur. Hence, $n$ can have neither square prime factors nor can be an even number(this is a Lehmer's theorem, which we prove it here by other method). This means that $t_{1}...t_{k}\leq 1$, thus some of $t_{1},...,t_{k}$ must be 0 or 1 or all of them be 1. Also, all of prime factors must be odd numbers. Certainly, if all of $t_{1},...,t_{k}$ be zero but one, then $n=p$ and the problem is solved. If the number of square-free prime factors are greater than or equal 2, then using Theorems 3 to 6 of this paper, we find out the equation
\begin{equation}
n-m\prod_{p_{i}\leq p\leq p_{k}}(p-1)=1
\end{equation}
does not any solution and Lehmer's totient conjecture is proven. 

\textbf{Acknowledgment}

This paper was submitted to a journal and quickly took some comments for its amendment. The author would like to thank anonymous referee for his/her nice comments. Also, thanks for Mr. Alexander Zujev, research scholar, University of California at Davis for his interesting comments on the paper. All the needed comments have been taken into account to the paper.


\begin{thebibliography}{99}
\bibitem{L} D.H. Lehmer, On Euler's totient function, Bull. Amer. Math. Soc.{38}(1932),745-751.
\bibitem{P} R. G. E. Pinch, A note on Lehmer's totient problem,personal web site, (2006) p.1.
\bibitem{Li} E. Lieuwens, Do there exist composite numbers for which $k\phi(M)=M-1$ holds?, Nieuw. Arch. Wisk.{18} (1970), 165-169. 
\bibitem{K} M. Kishore, On the number of distinct prime factors of $n$ for which $\phi(n)|(n-1)$,Nieuw. Arch. Wisk.{25} (1977), 48-53. 
\bibitem{CH} G.L. Cohen and P.Hagis jr., On the number of prime factors of $n$ if $\phi(n)|(n-1)$,Nieuw. Arch. Wiskd.,III. Ser. {28} (1980), 177-185. 
\bibitem{M}F. Mertens, Ein beitrag zur analytischen zahlentheorie, J.reine angew. Math. {78} (1874), 46-62.
\bibitem{WE} William and Fern Ellison,\textit{Prime Numbers},
Wiley-Interscience,a division of John Wiley and Sons.Inc.New York
(1985).
\bibitem{RS}J.B. Rosser, L.Schoenfeld,Approximate formulas for some functions of prime numbers, Illinois J.Math.\textbf{6} (1962) 64-94
\bibitem{Du} P. Dusart, Estimates of some functions over primes without R.H., arXiv:1002.0442v1, Feb. 2,2010
\end{thebibliography}
\end{document}